\documentclass[a4paper,notitlepage, 8pt,reqno]{article}

\usepackage[cp1251]{inputenc}
 \usepackage[T2A]{fontenc}
 \usepackage[russian,english]{babel}
\usepackage[all,arc,poly,2cell,curve,arrow,tips]{xy}
\usepackage{enumerate, eucal, amsthm,amsmath, amssymb}
\usepackage{graphicx}
\usepackage{mathtext}

\usepackage{tikz}

\allowdisplaybreaks[4]

\theoremstyle{definition}
\newtheorem{defn}{Definition}

\theoremstyle{plain}
\newtheorem{thm}{Theorem}
\newtheorem*{thm*}{Theorem}
\newtheorem{prop}{Proposition}
\newtheorem*{prop*}{Proposition}

\newtheorem{lem}{LEmma}

\newtheorem{remark}{В Remark}

\title{Classical  $6j$-symbols for finite dimensional representation of the algebra  $\mathfrak{gl}_3$ }

\author{D.V. Artamonov\\Lomonosov Moscow Stater Unversity\\artamonov@econ.msu.ru}

\begin{document}
	\maketitle
	
	\begin{abstract} 
		In the paper an explicit formula for an arbitrary   $6j$-symbol for finite-dimensional irreducible representations of the algebra  $\mathfrak{gl}_3$ is derived.  A $6j$-symbol is written as a result of substitution of   $\pm 1$ into a series of hypergeometric type, which is similar to a  $\Gamma$-series, which is a simplest example of a multivariate series of hypergeometric type.  The  selection rulers for a $6j$-symbol are derived.
		
	\end{abstract}
\section{Introduction}

\subsection{The Racah coefficients}

The Racah coefficients for a Lie algebra $g$ are defined as follows. Consider  irreducible representations $V^1$, $V^2$, $V^3$ of a given algebra and consider their tensor product $V^1\otimes V^2\otimes V^3$. Brackets in this work can be arranged in two ways:

$$
(V^1\otimes V^2)\otimes V^3 \text{ or } V^1\otimes (V^2\otimes V^3).
$$

One can split  $V^1\otimes V^2\otimes V^3$ into a sum of irreducible representations in two ways.

{\bf 1.} The first way. First split into irreducible $V^1\otimes V^2$:

\begin{equation}
\label{r1}
V^1\otimes V^2=\bigoplus_U Mult_{U}^{V^1,V^2}\otimes U,
\end{equation}
where  $U$ is an irreducible representation and  $Mult_{U}^{V^1,V^2}$ is a multiplicity space. This is a vector space which is not supplied with an action $g$. Take a tensor product of   \eqref{r1} and  $V^3$ from the right, one obtains

\begin{equation}
\label{rr1}
(V^1\otimes V^2)\otimes V^3=\bigoplus_{U,W} Mult_{U}^{V^1,V^2}\otimes Mult_{W}^{U,V^3}\otimes W
\end{equation}

{\bf 2.}  The second  way. First split into irreducible $V^2\otimes V^3$:

\begin{equation}
\label{r2}
V^2\otimes V^3=\bigoplus_U Mult_{H}^{V^2,V^3}\otimes H,
\end{equation}

and then

\begin{equation}
\label{rr2}
V^1\otimes (V^2\otimes V^3)=\bigoplus_{H,W} Mult_{H}^{V^2,V^3}\otimes Mult_{W}^{V^1,H}\otimes W
\end{equation}

There exists an isomorphism $\Phi: (V^1\otimes V^2)\otimes V^3 \rightarrow V^1\otimes (V^2\otimes V^3)$, which gives a mapping

\begin{equation}
\label{ph}
\Phi: \bigoplus_{U} Mult_{U}^{V^1,V^2}\otimes Mult_{W}^{U,V^3}\rightarrow  \bigoplus_{H} Mult_{H}^{V^2,V^3}\otimes
Mult_{W}^{V^1,H}
\end{equation}

\begin{defn}
 The Racah mapping is an induced by 
 $\Phi$ mapping 
	
	\begin{equation}
	\label{wrm}
	W\begin{Bmatrix}
	V^1 & V^2 & U\\ V^3 & W &H
	\end{Bmatrix}: Mult_{U}^{V^1,V^2}\otimes Mult_{W}^{U,V^3}\rightarrow  Mult_{H}^{V^2,V^3}\otimes Mult_{W}^{V^1,H}
	\end{equation}

\end{defn}

When one fixes  the basis, there appear matrix elements of this mapping. They are called {\it the Racah coefficients}. Let $s$ be an index listing the basis vectors in the multiplicity space/ Then one has a notation for the  Racah  coefficients

	\begin{equation}
\label{wrc}
W\begin{Bmatrix}
V^1 & V^2 & U\\ V^3 & W &H
\end{Bmatrix}^{s_1,s_2}_{s_3,s_4},
\end{equation}

where $s_1,s_2,s_3,s_4$ are indices of base vectors in  $Mult_{U}^{V^1,V^2}$, $Mult_{W}^{U,V^3}$, $  Mult_{Р }^{V^2,V^3}$,  $Mult_{W}^{V^1,H}$.

From this definition, the role of these coefficients from the point of view of representation theory is clear.
Let's consider a category of finite-dimensional representations and let us pass on to its Grothendieck ring. Equivalently, this means passing from the category of finite-dimensional  representations to the ring of their characters. Then some information about the category  is lost. For example,  one looses the information contained in Racah coefficients \cite{tk}. In some cases from  the Grothendieck ring and Racah coefficients, the category of representations can be restored.

The Racah  coefficients  have applications in quantum mechanics. They were introduced  by Racah in \cite{wgn}. In the case of $g=\mathfrak{sl}_2$, these coefficients are discussed in all textbooks on quantum mechanics, for example  \cite{wgn}, \cite{ll},  and textbooks on the theory of angular momentum, for example \cite{wrch}, \cite{bl}. In the case of finite-dimensional irreducible representations of $\mathfrak{sl}_2$, there are explicit formulas for the Racah coefficients in which they are expressed as the values of the hypergeometric function \cite{gk}. At the same time, the problem of calculating of  the Racah  coefficients for infinite-dimensional representations of $\mathfrak{sl}_2$ is still relevant \cite{d1}, \cite{d2}.

The Racah  coefficients for $\mathfrak{sl}_3$ are also important in applications
\cite{r81}, \cite{r86}.
The general formula for a $6j$-symbol for $\mathfrak{sl}_3$ is unknown (see  introduction to the dissertation \cite{slp}). Explicit formulas were obtained for some classes of representations \cite{r82}, \cite{r83}, \cite{r84}, \cite{r85}. In papers \cite{bl1}, \cite{bl2} there were also calculated some Racah coefficients for $\mathfrak{sl}_3$. These  Racah  coefficients play an important technical role in calculation of some Clebsch-Gordan coefficients for $\mathfrak{gl}_3$.

It should also be noted that even more attention is currently being paid to quantum Racah  coefficients, that is, the Racah  coefficients for quantum Lie algebras. It is not possible to give a complete overview of this activity, we  mention only recent works \cite{s1}, \cite{s2}, where the Racah  coefficients for algebras other than $U_q(\mathfrak{sl}_2)$ are considered.

 \subsection{$6j$-symbols}

These coefficients were introduced even earlier than the Racah coefficients by Wigner \footnote{these coefficients are not called the Wigner coefficients, since the term "Wigner coefficients" is used for the  coefficients that decompose a tensor product of two representations into irreducible ones, they are also Clebsch-Gordan coefficients}. To define them, we first introduce $3j$-symbols. Let $V^1$, $V^2$, $V^3$ be representations of $g$, $\{v^1_{\alpha}\}$, $\{v^2_{\beta}\}$, $\{v^3_{\gamma}\}$ are the bases in these representations.

  \begin{defn}  A $3j$-symbol is a collection of coefficients
  	
  	$$
  	\begin{pmatrix} V^1 & V^2 & V^3  \\ v^1_{\alpha} & v^2_{\beta}
  	& v^3_{\gamma} \end{pmatrix}^f\in \mathbb{C}, 
  	$$
  
 such that  $$\sum_{\alpha,\beta,\gamma} 	\begin{pmatrix} V^1 & V^2 & V^3  \\ v^1_{\alpha} & v^2_{\beta}
  & v^3_{\gamma} \end{pmatrix}^fv^1_{\alpha}\otimes v^2_{\beta}\otimes v^3_{\gamma}=f\in V^1\otimes V^2\otimes V^3,$$
  
 where  $f$ is a semi-invariant for the action of  $\mathfrak{gl}_3$. That is $f$ is an eigenvector for the action of the Cartan subalgebra and it vanishes under the action of  root elements.
  \end{defn}

  One can easily see a relation between $3j$-symbols and Clebsch-Gordan coefficients. To write it,  one takes a splitting \eqref{r1}. Let us  introduce a basis $e_s$ in the multiplicity space and put $U^s:=e_s\otimes U$, then \eqref{r1} one can be write  the following: $V^1\otimes V^2=\bigoplus_{s,U} U^s$. Choose the bases $\{v^1_{\alpha}\}$, $\{v^2_{\beta}\}$, $\{u^s_{\gamma}\}$ in these spaces. Then one  can write

 \begin{equation}
 \label{kg2}
 u_{\gamma}^s = \sum_{\alpha,\beta} D^{U,\gamma,s}_{V^1,V^2;\alpha,\beta}В v_{\alpha}^1\otimes v_{\beta}^2.
 \end{equation}

 The coefficient $D^{U,\gamma,s}_{V^1,V^2;\alpha,\beta}\in\mathbb{C}$ is called the Clebsch-Gordan coefficient, it is related to the $3j$-symbol by the rule

 $$
 D^{U,\gamma,s}_{V^1,V^2;\alpha,\beta}=\begin{pmatrix} V^1 & V^2 & \bar{U}  \\ v^1_{\alpha} & v^2_{\beta}
 & \bar{u}_{\gamma}  \end{pmatrix}^s,
 $$

where $\bar{U}$ and $\bar{u}_{\gamma}$ denote the contragradient representation and the dual basis in it. The index $s$  of a $3j$-symbol  carries the following information. The space of $3j$-symbols with given internal indexes is isomorphic to the multiplicity space. Therefore, by fixing the basis in the multiplicity space, one fixes a basic $3j$-symbols.

 \begin{defn}
 A $6j$-symbol is a paring of $3j$-symbols by the ruler:
 	
 	\begin{align}
 	\begin{split}
 	\label{6js}
& \begin{Bmatrix}
 V^1 & V^2 & U\\ V^3 & W &H
 \end{Bmatrix}^{s_1,s_2}_{s_3,s_4}:=\sum_{\alpha_1,...,\alpha_6}
 \begin{pmatrix} \bar{V}^1  & \bar{V}^2  & U \\ \bar{v}^1_{\alpha_1} & \bar{v}^2_{\alpha_2}  & u_{\alpha_4}\end{pmatrix}^{s_1}  \cdot \begin{pmatrix}  \bar{U}& \bar{V}^3  &W  \\\bar{ u}_{\alpha_4}  & \bar{v}^3_{\alpha_3} &w_{\alpha_5}   \end{pmatrix}^{s_2}\cdot \\ & \cdot  \begin{pmatrix}  V^2 &  V^3  &  \bar{H} \\ v^2_{\alpha_2} & v^3_{\alpha_3}   & \bar{h}_{\alpha_6}  \end{pmatrix}^{{s}_3}  \cdot \begin{pmatrix} V^1   & H  &\bar{W}  \\ v^1_{\alpha_1}   & h_{\alpha_6} &\bar{w}_{\alpha_5}\end{pmatrix}^{{s}_4}.
 	\end{split}
 	\end{align}
\end{defn}

This expression should be understood as follows: the Lie algebra $g$ acts  on the $3j$-symbols by acting on lower indices. We form a semiinvariant of $4$-x $3j$-symbols by pairing indexes so that in   two $3j$-symbols  only one pair of indexes paired.

Let's now write out a relationship between the Racah coefficients and $6j$-symbols. In the expression below, we use the fact that there is a duality between the spaces $Mult^{V^1,V^2}_U$ and $Mult^{\bar{V}^1,\bar{V}^2}_{\bar{U}}$. Thus if $s$ is an index of a basic vector in $Mult^{V^1,V^2}_U$ , then $\bar{s}$ is the index of the dual basic vector in $Mult^{\bar{V}^1,\bar{V}^2}_{\bar{U}}$.

$$
	W\begin{Bmatrix}
V^1 & V^2 & U\\ V^3 & W &H
\end{Bmatrix}^{\bar{s}_1,\bar{s}_2}_{s_3,s_4}=\begin{Bmatrix}
V^1 & V^2 & U\\ V^3 & W &H
\end{Bmatrix}^{s_1,s_2}_{s_3,s_4}
$$

Below we deal with   $6j$-symbols.

\subsection{The results of the paper }

In this paper, we obtain a simple explicit formula for an arbitrary $6j$ symbol for finite-dimensional irreducible representations of $\mathfrak{gl}_3$. This can be done by using the following ideas.

First, A-GKZ realization  of a representation of $\mathfrak{gl}_3$ is used. This realization is described in detail in \cite{smj}. The  space of a representation  is realized as a subspace in the space of polynomials in the variables $A_X$, $X\subset \{1,2,3\}$, antisymmetric in $X$, but not obeying other relations. The representation space is described as the space of polynomial solutions of a partial differential equation called the antisymmetrized Gelfand-Kapranov-Zelevinsky equation (A-GKZ for short). In this realization, it is possible to  write explicitly  functions in variables $A_X$ corresponding to the Gelfand-Tsetlin basis vectors (see the original work \cite{bb}, as well as \cite{smj}). It is important that in this model there is an explicitly written scalar product.

Secondly, an   explicitly construction  of  the multiplicity index   $s$ for the Clebsch-Gordan coefficients and $3j$-symbols  from  \cite{aa} is used.

Thirdly, an explicit and simple formula for an arbitrary $3j$-symbol obtained in \cite{smj} is use \footnote{in \cite{aa}, an explicit, but much more cumbersome solution of the same problem is obtained.}.

These ideas allow us to find an explicit simple formula for an arbitrary $6j$-symbol for finite-dimensional irreducible representations in terms of the values of a hypergeometric type function.

The structure of this paper is the following. The introductory Section \ref{r1} describes A-GKZ and functional realizations of representations of $\mathfrak{gl}_3$. A solution of the multiplicity problem in the decomposition of a the tensor product of $\mathfrak{gl}_3$ into the sum of irreducibles is given. A  a formula for a $3j$-symbol using  scalar products is presented.

In the main Section \ref{r30} an arbitrary $6j$-symbol for the algebra $\mathfrak{gl}_3$ is  calculated. First of all, in the Section \ref{r31} an implicit formula is given  via a scalar product (Lemma
\ref{ol1}). Then in the Section \ref{r32} the selection rules for the $6j$-symbol are proved (Theorem \ref{ot1}). Finally, in the section \ref{r33} an explicit formula is derived for the $6j$-symbol (Theorem \ref{ot2}).

\section{The basic notions}
\label{r1}
\subsection{  $A$-hypergeometric functions}

A detailed information about a $\Gamma$-series can be found in \cite{GG}.

Let $B\subset\mathbb{Z}^N$ be a lattice, let $\mu\in\mathbb{Z}^N$ be
a fixed vector. Define {\it a hypergeometric
$\Gamma$-series} in the variables $z_1,...,z_N$ by the formula

\begin{equation}
\label{gmr}
\mathcal{F}_{\mu}(z,B)=\sum_{b\in
	B}\frac{z^{b+\mu}}{\Gamma(b+\mu+1)},
\end{equation}
where $z=(z_1,...,z_N)$. We use a multi-index notation  in the numerator and denominator

$$
z^{b+\mu}:=\prod_{i=1}^N
z_i^{b_i+\mu_i},\,\,\,\Gamma(b+\mu+1):=\prod_{i=1}^N\Gamma(b_i+\mu_i+1).
$$

\begin{remark}
A more compact formula for  \eqref{gmr} is the following: $\mathcal{F}_{\mu}(z,B)=\sum_{x\in\mu+B}\frac{z^{x}}{\Gamma(x+1)}$.
\end{remark}

Note that if at least one of the components of the vector $b+\mu$ is a negative integer, then the corresponding term in \eqref{gmr} turns to zero. Due to this, there is be only a finite number of terms in the series considered in the paper. For simplicity, we write factorials instead of $\Gamma$-functions.

A $\Gamma$-series satisfies the Gelfand-Kapranov-Zelevinsky system (the GKZ system for short). Let us write it explicitly in the case  $z=(z_3,z_1,z_2,z_{1,3},z_{2,3},z_{1,2})$, $B=\mathbb{Z}<v=e_1-e_2-e_{1,3}+e_{2,3}>$:
\begin{align}
\begin{split}
\label{gkzs} \Big(\frac{\partial^2}{\partial z_{1}\partial
	z_{2,3}}-\frac{\partial^2}{\partial z_2\partial z_{1,3}}\Big)\mathcal{F}_{\mu}(z,B)&=0,
\\
(z_1\frac{\partial}{\partial z_1}+ z_2\frac{\partial}{\partial
	z_2})\mathcal{F}_{\mu}(z,B)&=(\mu_1+\mu_2)\mathcal{F}_{\mu}(z,B),\quad
\\
(z_1\frac{\partial}{\partial z_1}+
z_{1,3}\frac{\partial}{\partial z_{1,3}})\mathcal{F}_{\mu}(z,B)&=(\mu_1+\mu_{1,3})\mathcal{F}_{\mu}(z,B),
\\
(z_1\frac{\partial}{\partial z_1}-z_{2,3}\frac{\partial}{\partial
	z_{2,3}})\mathcal{F}_{\mu}(z,B)&=(\mu_1-\mu_{2,3})\mathcal{F}_{\mu}(z,B),\\
z_3\frac{\partial}{\partial z_3}\mathcal{F}_{\mu}(z,B)=\mu_{3}\mathcal{F}_{\mu}(z,B),\,\,\,\,
&\,\,\,\,\,  z_{1,2}\frac{\partial}{\partial z_{1,2}}\mathcal{F}_{\mu}(z,B)=\mu_{1,2}\mathcal{F}_{\mu}(z,B).
\end{split}
\end{align}

\subsection{ A-GKZ realization of a representation of   $\mathfrak{gl}_3$}

Throughout the rest of the work, we  talk about finite-dimensional irreducible representations.

Details can be found in \cite{smj}. Consider the variables $A_X$, where $X\subset \{1,2,3\}$ is a proper subset, antisymmetric in $X$, but not obeying other relations. On these variables there is an action $\mathfrak{gl}_3$, defined by the rule

\begin{equation}
\label{edeta}
E_{i,j}A_X=\begin{cases}  A_{X\mid_{j\mapsto i}},\,\,\,\,j\in X\\ 0\text{ otherwise }\end{cases}
\end{equation}
 
 Here $X\mid_{j\mapsto i}$ denotes a substitution of  $j$ by  $i$.
 
 Consider an equation called the antisymmetrized GKZ equation (the A-GKZ equation):
 
 \begin{equation}
 \label{agkz}
(\frac{\partial^2}{\partial A_{1}A_{2,3}}-\frac{\partial^2}{\partial A_{2}A_{1,3}}+\frac{\partial^2}{\partial A_{3}A_{1,2}})F=0
 \end{equation}
 
 \begin{remark}

 The term A-GKZ is explained as follows. Consider the system \eqref{gkzs}, leave only the first equation   and "antisymmetrize" it by adding a third term. The resulting equation is the equation A-GKZ.
 \end{remark}
 
 One has.
 
 \begin{thm}[\cite{smj}]
 	\label{t1}
The space of polynomial solutions \eqref{agkz} is invariant under the action of $\mathfrak{gl}_3$. As a representation, the space of polynomial solutions is a direct sum with multiplicity $1$ of all finite-dimensional irreducible representations  with $m_{3}=0$. The space of an irreducible representation with the highest weight $[m_{1}, m_{2},0]$ has the highest vector $$A_1^{m_1-m_2}A_{1,2}^{m_2}$$ and consists of all polynomial solutions whose homogeneous degree in variables $A_X$, $|X|=1$ is equal to $m_{1}-m_2$, and homogeneous degree in variables $A_X$, $|X|=2$ is equal to $m_2$.
 \end{thm}
 
Let us write a base in the space of polynomial solutions. Consider the space    $\mathbb{C}^6$ with  coordinates  $A_1,A_2,A_3,A_{1,2},A_{1,3},A_{2,3}$ and define the vectors 
 
 \begin{align}
 \begin{split}
&v=e_1-e_2-e_{1,3}+e_{2,3},\,\,\,r=e_3+e_{1,2}-e_1-e_{2,3},
\end{split}
\end{align}
Now consider a  $\Gamma$-series associated with a lattice  $B=\mathbb{Z}<v>$ and a vector  $\mu\in \mathbb{Z}^6$:

\begin{equation}
\mathcal{F}_{\mu}(A,B):=\sum_{t\in\mathbb{Z}}\frac{A^{\mu+tv}}{(\mu+tv)!},
\end{equation}

where we used the multi-index notations $$A^{\mu+tv}=\prod_X A_X^{\mu_X+tv_X},\,\,\,\,(\ mu+tv)!=\prod_X (\mu_X+tv_X)!.$$ Then (see calculations in \cite{aa}, \cite{a5} or \cite{smj}) the basis in the space of polynomial solutions consists of non-zero functions of type

  \begin{equation}
 F_{\mu}(A):=\sum_{s\in\mathbb{Z}_{\geq 0}}  q^{\mu}_s \zeta_A^{s}\mathcal{F}_{\mu-sr}(A), 
 \end{equation}
 
 where  \begin{align}\begin{split}\label{cs} &t_0=1,В \,\,\,\,\, t^{\mu}_s=\frac{1}{s(s+1)+s(\mu_1+\mu_2+\mu_{1,3}+\mu_{2,3})}, \text{ in the case }s>0\\& q_{\mu}^s=\frac{t_s^{\mu}}{\sum_{s'\in\mathbb{Z}_{\geq 0}}  t_{s'}^{\mu}}  \\&r=e_3+e_{1,2}-e_1-e_{2,3},\,\,\,\,\,
 \zeta_A=A_1A_{2,3}-A_{2}A_{1,3}.
 \end{split}\end{align}

Note that actually $F_{\mu}(A)$   depends not on   the vector  $\mu$, but on the shifted lattice   $\Pi=\mu+\mathbb{Z}<v>$.

 Moreover, when this function is nonzero, then it is nothing else but the Gelfand-Tsetlin basic vector. Namely, the Gelfand-Cetlin diagram $$(m_{p,q})=\begin{pmatrix} m_{1,3 } &&m_{2,3} &&m_{3,3} \\&m_{1,2} &&m_{2,2} \\&&m_{1,1}\end{pmatrix}{}$$ corresponds to such a function $F_{\mu}(A)$, that the shifted lattice $\Pi=\mu+\mathbb{Z}<v>$ is given by the equations

 \begin{equation}
 \delta\in \Pi \Leftrightarrow\begin{cases}
 \delta_1+\delta_2+\delta_3+\delta_{1,2}+\delta_{1,3}+\delta_{2,3}=m_{1,3},\\
 \delta_{1,2}+\delta_{1,3}+\delta_{2,3}=m_{2,3},\\
  \delta_1+\delta_2+\delta_{1,2}+\delta_{1,3}+\delta_{2,3}=m_{1,2},\\
   \delta_{1,2}=m_{2,2},\\
    \delta_1+\delta_{1,2}+\delta_{1,3}=m_{1,1}
 \end{cases}
 \end{equation}

Moreover, this is a one-to-one correspondence between  the Gelfand-Tsetlin diagrams and nonzero functions $F_{\mu}(A)$.

It is important that in the A-GKZ realization there is an explicitly formula for an invariant scalar product

 \begin{equation}
 \label{skp}
 <f(A),g(A)>=f(\frac{\partial }{\partial A})g(A)\mid_{A=0}.
 \end{equation}

 Also note that if the representation of $V$ is realized  in the space of polynomials in the variables $A_X$: $V=\{h(A)\}$, then the contragradient representation is realized in the space of polynomials in the operators $\frac{\partial}{\partial A_X}$.  The action  of$\mathfrak{gl}_3$ on the differential operators is generated by an the action on functions of $A_X$. In this case $\bar{V}=\{h(\frac{\partial }{ \partial A_X}): \,\,\, h(A)\in V\}$. The pairing is given by a formula similar to \eqref{skp}:

 \begin{equation}
 \label{spar}
  \{h_1(A),h_2(\frac{\partial }{ \partial A_X})\}>=h_2(\frac{\partial }{ \partial A_X})h_1(A)\mid_{A=0}.
 \end{equation}

 Since the base  $F_{\mu}(A)$ is the Gelfand-Tsetlin base, which is orthogonal, then the dualto $F_{\mu}(A)$   base  is the base $\frac{1}{|F_{\mu}|^2}F_{\mu}(\frac{\partial }{\partial A})$.

 \subsection{The functional realization}
 \label{fr}
 We need also another realization, details about it can be found in \cite{zh}.
Functions on the group $GL_3$ form a representation of $GL_3$. On
a function $f(g)$, $g\in GL_3$, an element of the group $X\in GL_{3}$ acts
by right shifts according to the rule
 
 \begin{equation}
 \label{xf} (Xf)(g)=f(gX).
 \end{equation}

 Passing to an  infinitesimal version of this action,  one gets that the space
of all functions on $GL_3$  is a representation of $\mathfrak{gl}_3$.

Any finite-dimensional irreducible representation can be
embedded into  the  space of function on $\mathfrak{gl}_3$.
More precise, if
$[m_{1},m_2,m_{3}]$ is the highest weight, then in the space of all functions
there is a highest vector with such a weight, which is explicitly written as follows.

Let $a_{i}^{j}$, $i,j=1,2,3$ be
a function of a matrix element on the group $GL_{3}$. Here $j$ is the row index and $i$ is the column index.
Also put

 \begin{equation}
 \label{dete}
 a_{i_1,...,i_k}:=det(a_i^j)_{i=i_1,...,i_k}^{j=1,...,k},
 \end{equation}

where  one takes a determinant  of the submatrix in the matrix $(a_i^j)$
formed by rows indexed by  $1,...,k$ and columns
$i_1,...,i_k$.
The operator $E_{i,j}$ acts on a determinant by
acting onto column indexes, according to a formula similar to \eqref{edeta}

\begin{equation}
\label{edet1}
E_{i,j}a_X=\begin{cases}  a_{X\mid_{j\mapsto i}},j\in X\\ 0\text{ otherwise.}\end{cases}
\end{equation}

There is a mapping from A-GKZ realization to a functional one, consisting in replacing \begin{equation}\label{post}A_X\mapsto a_x.\end{equation}

The fundamental difference between this realization and the A-GKZ realizattion is that the determinants of $a_x$ obey the Plucker relation

\begin{equation}
\label{spl}
a_{1}a_{2,3}-a_2a_{1,3}+a_3a_{1,2}=0,
\end{equation}

and the variables   $A_X$ are independent.

On the one hand, this simplifies the description of the representation space. So, there is a Theorem (cf. with the Theorem \ref{t1}).

\begin{thm}
The space of an irreducible representation with the highest weight $[m_{1}, m_{2},0]$ and the highest vector $$a_1^{m_1-m_2}a_{1,2}^{m_2}$$  consists of all {\it polynomials} whose homogeneous degree in variables $a_X$, $|X|=1$ is equal to $m_{1}-m_2$, and the homogeneous degree in variables $a_X$, $|X|=2$ is equal to $m_2$.
\end{thm}

But on the other hand, many calculations become more complicated.
So the formula for the invariant scalar product \eqref{skp} is incorrect when $A_X$ is mechanically replaced by $a_X$.

Nevertheless, this following is true. When substituting \eqref{post}, the function $F_{\mu}(A)$ passes into $\mathcal{F}_{\mu}(a,B)$. So $\mathcal{F}_{\mu}(a,B)$ is a Gelfand-Tsetlin basic vector in a functional realization.

\subsection{ A solution of the multiplicity problem  for  $3j$-symbols }
\label{3jkr}

Let us construct  explicitly  the multiplicity index of $s$ for $3j$-symbols.
The triple tensor product can be realized  in the function space on $GL_3\times GL_3\times GL_3$. Functions of  matrix elements on these factors $GL_3$ are denoted as $a_i^j$, $b_i^j$, $c_i^j$. Similar letters will denote the determinants of matrices composed of these matrix elements.

Let us construct explicitly a semiinvariant vectors $f$. It is easy to understand that  semiinvariants are the following expressions

 \begin{equation}
\label{skob}
f=\frac{(abc)^{\tau_1}(aac)^{\tau_2}(acc)^{\tau_3}(aab)^{\tau_4}
(abb)^{\tau_5}(bbc)^{\tau_6}(bcc)^{\tau_7}(aabbcc)^{\tau_8}}{\tau_1!\tau_2! \tau_3!\tau_4!\tau_5!\tau_6!\tau_7!\tau_8!},\end{equation}

where 
\begin{equation}
(abc)=det\begin{pmatrix}
a_1^1 & a_2^1 & a_3^1\\
b_1^1 & b_2^1 & b_3^1\\
c_1^1 & c_2^1 & c_3^1
\end{pmatrix},\,\,\, (aac)=det\begin{pmatrix}
a_1^1 & a_2^1 & a_3^1\\
a_1^2 & a_2^2 & a_3^2\\
c_1^1 & c_2^1 & c_3^1
\end{pmatrix}\text{ Рё.С‚.Рґ.}
\end{equation}

Here \begin{align*}&(aabbcc):=(\tilde{a}\tilde{b}\tilde{c}),\,\,\,\tilde{a}_1^1:=a_{2,3},\,\,\,\tilde{a}_2^1:=-a_{1,3},\,\,\,\tilde{a}_3^1:=a_{1,2},\end{align*}
$\tilde{b}_i^1$,  $\tilde{c}_i^1$ are defined analogously.

In this case, $f$ belongs to  the tensor product  $V^1\otimes V^2\otimes V^3$  of representations with highest weights $[m_1,m_2,0]$, $[m'_1,m'_2,0]$, $[M_1,M_2,0]$ if and only if the  following conditions hold

\begin{align}
\begin{split}\label{pvl}
m_1=\tau_1+\tau_2+\tau_3+\tau_4+\tau_5+\tau_8,\,\,\, m_2=\tau_2+\tau_4+\tau_8,\\
m'_1=\tau_1+\tau_4+\tau_5+\tau_6+\tau_7+\tau_8,\,\,\, m'_2=\tau_5+\tau_6+\tau_8,\\
M_1=\tau_1+\tau_2+\tau_3+\tau_6+\tau_7+\tau_8,\,\,\,M_2=\tau_3+\tau_7+\tau_8
\end{split}
\end{align}

\begin{prop}[see Proposition 2 in  \cite{aa}]
	\label{mlt}
In the space of  $3j$-symbols with the same internal indexes, there is a set of generators  consisting of $3j$-symbols indexed by semiinvariant functions $f$ of the form \eqref{skob}, consistent with the higher weights in the upper row of the $3j$-symbol according to the rule \eqref{pvl}. To get a basis, it is necessary to leave $f$ such that either $\tau_1=0$ or $\tau_8=0$.
\end{prop}

\subsection{A formulas for a  $3j$-symbol}

In  representations, we take bases of the type $F_{\mu}(A)$, but in order to coordinate the notation with \eqref{6js}, we  use the index $\alpha_i$ to enumerate the basis vectors. Thus, as $v^1_{\alpha_1}$ we  take $F_{\alpha_1}(A^1)$, as $v^2_{\alpha_2}$ we take $F_{\alpha_2}(A^2)$, etc.

Also, as a multiplicity index for the $3j$-symbol, we  use a function defines a semiinvariant.

The calculations in the present paper are based on the formula for a $3j$-symbols in the basis  $F_{\alpha}$ obtained in \cite{smj}.
Let $G$ be a semiinvariant in $V^1\otimes V^2\otimes V^3$ in the A-GKZ realization  and $G(A,B,C)=\sum_{\alpha'_1,\alpha'_2,\alpha'_3} c_{\alpha'_1,\alpha'_2,\alpha'_3} F_{\alpha'_1}(A)\otimes F_{\alpha'_2}(B)\otimes F_{\alpha'_3}(C)$, where

$$
c_{\alpha_1,\alpha_2,\alpha_3}=\begin{pmatrix}
V^1 & V^2 & V^3\\ F_{\alpha_1} &  F_{\alpha_2} & F_{\alpha_3}
\end{pmatrix}^G
$$
The vectors $F_{\alpha}$ are orthogonal, since they are vectors of the Gelfand-Tsetlin basis, so

$$
c_{\alpha_1,\alpha_2,\alpha_3}=\frac{<G(A,B,C),F_{\alpha_1}(A) F_{\alpha_2}(B) F_{\alpha_3}(C)>}{|F_{\alpha_1}(A)|^2| F_{\alpha_2}(B)|^2 |F_{\alpha_3}(C)|^2}
$$

Also note the following. Let the semiinvariant $G$ in the A-GKZ  realization correspond to a semiinvariant $f$ in the functional realization given by the function \eqref{skob}. Then

 $$
 G(A,B,C)=f(A,B,C)+r(A,B,C),
 $$ 

where $r(A,B,C)$  belongs to the ideal of  relations between the determinants $a_X$, $b_X$, $c_X$. This ideal is generated by the relations \eqref{spl}, as well as similar relations for $b_X$, $c_X$. Then

$$
<r(A,B,C),F_{\alpha_1}(A) F_{\alpha_2}(B) F_{\alpha_3}(C)>=0,
$$
since $F_{\alpha_1}(A)$, $F_{\alpha_2}(B)$, $F_{\alpha_3}(C)$ are the solution of A-GKZ equations
\eqref{agkz}, and the scalar product is given by \eqref{skp}. So, a $3j$-symbol is calculated in the following manner

 \begin{equation}
 \label{s3j}
 \begin{pmatrix}
 V^1 & V^2 & V^3\\ F_{\alpha_1} &  F_{\alpha_2} & F_{\alpha_3}
 \end{pmatrix}^f=\frac{<f(A,B,C),F_{\alpha_1}(A) F_{\alpha_2}(B) F_{\alpha_3}(C)>}{|F_{\alpha_1}(A)|^2| F_{\alpha_2}(B)|^2 |F_{\alpha_3}(C)|^2}
 \end{equation}

\section{ $6j$-symbols}

\label{r30}

Let's start calculating  the $6j$-symbol \eqref{6js}. First, we  give an expression for it in terms of scalar products, and then we  give an explicit formula through the value of a hypergeometric type function, with $\pm 1$ substituted as arguments.

 \subsection{ An expression though scalar products }
 
 \label{r31}

 Let's write an expression for the $6j$ symbol \eqref{6js}. Let the multiplicity indices $s_1,...,s_4$ in \eqref{6js} correspond to the functions $f_1,...,f_4$.

We need to calculate  $3j$-symbols for the contragradient representation and the dual basis. Let use a realization of a contragradient representation described at the end of the section \ref{fr}. The a basis dual to $F_{\alpha_i}(A^i)$ is the basis of $\frac{1}{|F_{\alpha_i}|^2}F_{\alpha_i}(\frac{\partial}{\partial A^i})$.

Note that a $3j$-symbol of the form  $$\begin{pmatrix} \bar{V}^1 & \bar{V}^2 & U\\F_{\alpha_1}(\frac{\partial}{\partial A^1}) &  F_{\alpha_2}(\frac{\partial}{\partial A^2}) & F_{\alpha_4}(A^4)   \end{pmatrix}^f$$
can be calculated as follows:

\begin{equation}
\label{s3j1}
\begin{pmatrix} \bar{V}^1 & \bar{V}^2 & U\\F_{\alpha_1}(\frac{\partial}{\partial A^1}) &  F_{\alpha_2}(\frac{\partial}{\partial A^2}) & F_{\alpha_4}(A^4)   \end{pmatrix}^f=\frac{<f(\frac{\partial}{\partial A^1},\frac{\partial}{\partial A^2},A^4),F_{\alpha_1}(\frac{\partial}{\partial A^1}) F_{\alpha_2}(\frac{\partial}{\partial A^2}) F_{\alpha_4}(A^4)>}{|F_{\alpha_1}(\frac{\partial}{\partial A^1})|^2| F_{\alpha_2}(\frac{\partial}{\frac{\partial}{\partial A^2}})|^2 |F_{\alpha_4}(A^4)|^2}.
\end{equation}

The scalar product in case when there is a function not of  a variable, but of a differentiation operator is calculated using a formula similar to \eqref{skp}. One has

\begin{align}
\begin{split}
\label{ss1}
&<f(\frac{\partial}{\partial A^1},\frac{\partial}{\partial A^2},A^4),F_{\alpha_1}(\frac{\partial}{\partial A^1}) F_{\alpha_2}(\frac{\partial}{\partial A^2}) F_{\alpha_4}(A^4)>=
	<f( A^1, A^2,A^4),F_{\alpha_1}(A^1) F_{\alpha_2}(\partial A^2) F_{\alpha_4}(A^4)>,\\
	& |F_{\alpha_1}(A^1)|^2=|F_{\alpha_1}(\frac{\partial}{\partial A^1})|^2,...
	\end{split}
\end{align}

Bases $F_{\alpha_1}(A^1)$ and $F_{\alpha_1}(\frac{\partial}{\partial A^1})$ etc. are not dual,  the basis  dual to $F_{\alpha_1}(A^1)$ is  $\frac{1}{|F_{\alpha_1}|^2}F_{\alpha_1}(\frac{\partial}{\partial A^1})$. So the $6j$-symbol is expressed in terms of the considered $3j$-symbols \eqref{s3j1} as follows

\begin{align}
\begin{split}
\label{6js1}
& \begin{Bmatrix}
V^1 & V^2 & U\\ V^3 & W &H
\end{Bmatrix}^{f_1,f_2}_{f_3,f_4}:=\sum_{\alpha_1,...,\alpha_6}
\begin{pmatrix} \bar{V}^1 & \bar{V}^2 & U\\F_{\alpha_1}(\frac{\partial}{\partial A^1}) &  F_{\alpha_2}(\frac{\partial}{\partial A^2}) & F_{\alpha_4}(A^4)   \end{pmatrix}^{f_{1}}  \cdot\\&\cdot \begin{pmatrix}  \bar{U}& \bar{V}^3  &W  \\   F_{\alpha_4}(\frac{\partial}{\partial A^4})  &   F_{\alpha_3}(\frac{\partial}{\partial A^3}) & F_{\alpha_4}(A^5) \end{pmatrix}^{f_2}\cdot \\ & \cdot  \begin{pmatrix}  V^2 &  V^3  &  \bar{H} \\  F_{\alpha_2}(A^2) & F_{\alpha_3}(A^3)   & F_{\alpha_6}(\frac{\partial}{\partial A^6}) \end{pmatrix}^{{f}_{3}}  \cdot \begin{pmatrix} V^1   & H  &\bar{W}  \\  F_{\alpha_1}(A^1)   &  F_{\alpha_6}(A^6)  &  F_{\alpha_5}(\frac{\partial}{\partial A^5})\end{pmatrix}^{{f}_{4}} \cdot |F_{\alpha_1}|^2\cdot....\cdot |F_{\alpha_6}|^2  
\end{split}
\end{align}

Take the expressions \eqref{s3j1} and substitute them in \eqref{6js1}. Consider \eqref{ss1}. At the same time the expression $|F_{\alpha_i}|^2$  occurring 
at the end of \eqref{6js}   are written as $F_{\alpha_i}(\frac{\partial}{\partial A^i})F_{\alpha_i}(A^i)\mid_{A=0}$
В 

One obtains В  

В  
В  \begin{align*}
&В   \begin{Bmatrix}
В  V^1 & V^2 & U\\ V^3 & W &H
В  \end{Bmatrix}^{f_1,f_2}_{f_3,f_4}=\sum_{\alpha_1,...,\alpha_6}\frac{<f_{1},F_{\alpha_1}F_{\alpha_2}F_{\alpha_4}>}{  |F_{\alpha_1}|^2     |F_{\alpha_2}|^2   |F_{\alpha_4}|^2 }   F_{\alpha_1}(\frac{\partial}{\partial A^1})F_{\alpha_2}(\frac{\partial}{ \partial A^2})F_{\alpha^4}(A^4) \cdot\\& \frac{<f_{2},F_{\alpha_4}F_{\alpha_3}F_{\alpha_5}>}{  |F_{\alpha_4}|^2     |F_{\alpha_3}|^2   |F_{\alpha_5}|^2 }   F_{\alpha_4}(\frac{\partial}{\partial A^4})F_{\alpha_3}(\frac{\partial}{ \partial A^3})F_{\alpha_5}(A^5)...\mid_{A_1=...=A_6=0}      
В  \end{align*}

Now write
  
$$ f_1(  \frac{\partial}{\partial A^1}, \frac{\partial}{\partial A^2} , A^4 )=\sum  \frac{<f_{1},F_{\alpha_1}F_{\alpha_2}F_{\alpha_4}>}{  |F_{\alpha_1}|^2     |F_{\alpha_2}|^2   |F_{\alpha_4}|^2 }   F_{\alpha_1}(\frac{\partial}{\partial A^1})F_{\alpha_2}(\frac{\partial}{ \partial A^2})F_{\alpha^4}(A^4),
$$
and analogous expressions forВ  В $ f_2(  \frac{\partial}{\partial A^4}, \frac{\partial}{\partial A^5} , A^5 )$ , $f_3(  \frac{\partial}{\partial A^2},  A^3, \frac{\partial}{\partial A^6}  )$,   $f_4(  A^1,  A^6, \frac{\partial}{\partial A^5}  )$.
 
 Using thatВ  $\{F_{\alpha_i}(A^i),   F_{\alpha'_i}(\frac{\partial}{\partial A^i}) \}=|F_{\alpha_i}|^2$ if   $\alpha_i=\alpha'_i$ and  $0$ otherwise one gets

\begin{lem}
	\label{ol1}
	\begin{align}
	\begin{split}
	\label{6jsk}
	&\begin{Bmatrix}
	V^1 & V^2 & U\\ V^3 & W &H
	\end{Bmatrix}^{f_1,f_2}_{f_3,f_4}= f_1(  \frac{\partial}{\partial A^1}, \frac{\partial}{\partial A^2} , A^4 )f_2(  \frac{\partial}{\partial A^4}, \frac{\partial}{\partial A^3} , A^5 )\cdot\\&\cdot f_3(  \frac{\partial}{\partial A^2},  A^3, \frac{\partial}{\partial A^6}  )f_4(  A^1,  A^6, \frac{\partial}{\partial A^5}  ).\mid_{A_1=...=A_6=0}  ,
	\end{split}
	\end{align}

 where it is assumed that this expression is evaluated as follows. At first it is represented as the sum of monomials of variables and differential operators. We multiply variables and differential operators as if they commute. Then, in each monomial, the differential operators act on the product of the variables occurring in the monomial. Finally, zeros are substituted instead of all the variables.
\end{lem}

\subsection{ The selection rulers}

\label{r32}

Let's find the conditions necessary for the $6j$ symbol to be nonzero. To do this, we consider supports of various functions.

\begin{defn}

A support  of a function $f$ of variables $Z=\{z_1,...,z_N\}$ is defined as follows. Decompose the function into a power series, take the set of exponents $\delta\in\mathbb{Z}^N$ of  monomials $Z^{\delta}:=z_1^{\delta_1}\cdot...\cdot z_N^{\delta_N}$, occurring in the decomposition of $f$. The obtained set of exponents is called the support of the function.

The support is denoted as  $supp f$.
\end{defn}

In \cite{smj} the function $f$ of type \eqref{skob} is considered as a function of  ther variables

\begin{align}
\begin{split}
\label{pz}
&
Z=\{[c_1a_{2,3}] ,[c_2a_{1,3}] ,[c_3a_{1,2}] , [a_1c_{2,3}] , [a_{2}c_{1,3}] ,[a_3c_{1,2}] ,
[c_1b_{2,3}]  ,  [c_2b_{1,3}] ,\\& [c_3b_{1,2}] , [b_1c_{2,3}] ,[b_{2}c_{1,3}] ,[b_3c_{1,2}] ,
[b_1a_{2,3}]  ,  [b_2a_{1,3}] , [b_3a_{1,2}], [a_1b_{2,3}] , [a_{2}b_{1,3}] ,\\& [a_3b_{1,2}] ,
[a_1b_2c_3],[a_2b_3c_1], [a_3b_1c_2], -[a_2b_1c_3],[a_1b_3c_2], [a_3b_2c_1],\\&
[a_{2,3}b_{1,3}c_{1,2}]  ,    [a_{1,3}b_{1,2}c_{2,3}]  ,    [a_{1,2}b_{2,3}c_{1,3}]  , 
-[a_{1,3}b_{2,3}c_{1,2}]  ,   [a_{2,3}b_{1,2}c_{1,3}]  ,   [a_{1,2}b_{1,3}c_{2,3}] \}. 
\end{split}
\end{align}

Thus the  summands in the determinants $(cas),(ac),...,(aabbcc)$ are identified with these variables.

In  \cite{smj} it is shown that the support of the function $f$ of type \eqref{skob} has the form of an intersection of a shifted lattice and a positive octant

$$
suppf=\Big(   \bigcap_i(\delta_i\geq 0)  \Big) \cap (\kappa+\mathcal{B})\subset \mathbb{Z}^{30},
$$

where  $30$  is the number of variables   $Z$.  The lattice $\mathcal{B}\subset\mathbb{Z}^{30}$ is generated by the vectors

\begin{align*}
&p_1= e_{[c_1a_{2,3}]}-e_{ [c_2a_{1,3}] }, \,\,\, \,\,\,  p_2=e_{[c_1a_{2,3}]}-e_{ [c_3a_{1,2}] },\,\,\, \,\,\,     p_3= e_{[a_1c_{2,3}] }-e_{[a_2c_{1,3}] }, \\& p_4= e_{[a_1c_{2,3}] }-e_{[a_3c_{1,2}] }, \,\,\, \,\,\, 
p_5= e_{[c_1b_{2,3}]}-e_{ [c_2b_{1,3}] },  \,\,\, \,\,\,   p_6=e_{[c_1b_{2,3}]}-e_{ [c_3b_{1,2}] },\\&p_7= e_{[b_1c_{2,3}] }-e_{[b_2c_{1,3}] },   \,\,\, \,\,\,  p_8= e_{[b_1c_{2,3}] }-e_{[b_3c_{1,2}] }, \,\,\, \,\,\,    p_9= e_{[a_1b_{2,3}]}-e_{ [a_2b_{1,3}] },\\& p_{10}=e_{[a_1b_{2,3}]}-e_{ [a_3b_{1,2}] }, \,\,\, \,\,\,   p_{11}= e_{[b_1a_{2,3}] }-e_{[b_2a_{1,3}] }, \,\,\, \,\,\,   p_{12}= e_{[b_1a_{2,3}] }-e_{[b_3a_{1,2}] },\\&p_{13}=e_{[a_1b_2c_3]}-e_{[a_2b_3c_1]} , \,\,\, \,\,\,     p_{14}=e_{[a_1b_2c_3]}-e_{[a_3b_1c_2]},   \,\,\, \,\,\, 
p_{15}=e_{[a_1b_2c_3]}-e_{[a_2b_1c_3]} \\&p_{16}=e_{[a_1b_2c_3]}-e_{[a_1b_3c_2]},  \,\,\, \,\,\,  p_{17}=e_{[a_1b_2c_3]}-e_{[a_3b_2c_1]},   \,\,\, \,\,\, 
p_{18}= e_{[a_{2,3}b_{1,3}c_{1,2}]}-e_{ [a_{1,3}b_{1,2}c_{2,3}] }, \\&   p_{19}= e_{[a_{2,3}b_{1,3}c_{1,2}]}-e_{ [a_{1,2}b_{2,3}c_{1,3}] },  \,\,\, \,\,\,    p_{20}= e_{[a_{2,3}b_{1,3}c_{1,2}]}-e_{ [a_{1,3}b_{2,3}c_{1,2}] },   \,\,\, \,\,\,  \\&  p_{21}= e_{[a_{2,3}b_{1,3}c_{1,2}]}-e_{ [a_{1,2}b_{1,3}c_{2,3}] },   \,\,\, \,\,\, p_{22}= e_{[a_{2,3}b_{1,3}c_{1,2}]}-e_{ [a_{2,3}b_{1,2}c_{1,3}] },
\end{align*}

 and  $\kappa$ is the vector
 
 $$
\kappa=(\tau_2,0,0,\tau_3,0,0,\tau_6 ,0,0,\tau_7,0,0,\tau_4,0,0,\tau_5,0,0,\tau_1,0,0,0,0,0,\tau_8,0,0,0,0,0)
 $$

Moreover, $f$ is $\Gamma$ series in variables $Z$, taken with the sign $\pm 1$, constructed from the lattice $\mathcal{B}$ and the shift vector $\kappa$. One takes a variable with the $-$ sign if it occurs in a determinant with the minus sign.

The function $f$ can also be understood in a more familiar way  as a function of $A_X$, $B_X$, $C_X$. The transition from variables $Z$ to variables $A_X$, $B_X$, $C_X$ is carried out using an obvious substitution

\begin{equation}
\label{zmn}
[c_1a_{2,3}]\mapsto C_1A_{2,3},...
\end{equation}

There appears mappings $pr_a, pr_b, pr_c$  from the space of exponents for variables $Z$ to the space of exponents for variables $A_X$, $B_X$, $C_X$, respectively:

$$
pr_a,pr_b,pr_c:\mathbb{Z}^{30}\rightarrow \mathbb{Z}^6.
$$

Here $\mathbb{Z}^{30}$ is a lattice of exponents of monomials in variables $Z$, and $\mathbb{Z}^6$ is a lattice of exponents of monomials in variables $A_X$ (or $B_X$, or $C_X$). In  $\mathbb{Z}^{30}$ one has a basis indexed by the variables $Z$, and $\mathbb{Z}^6$ has a basis indexed by proper subsets of $X\subset\{1,2,3\}$.

Now we can formulate the condition necessary for the $6j$-symbol to be nonzero. Consider the expression \eqref{6jsk}. Let $$supp f_i=\Big ( \bigcap_i(\delta_i\geq 0) \Big ) \cap (\kappa_i+\mathcal{B}_i)\subset\mathbb{Z}^{30}, i=1,...,4,$$ put

\begin{equation}
\label{h}
H= suppf_1\oplus supp f_2\oplus   suppf_3\oplus supp f_4 \subset \mathbb{Z}^{30} \oplus   \mathbb{Z}^{30} \oplus   \mathbb{Z}^{30}\oplus   \mathbb{Z}^{30}
\end{equation}

For each $i=1,...,4$ there is a mapping $pr^i_a,pr^i_b,pr^i_c : \mathbb{Z}^{30}\rightarrow \mathbb{Z}^{6}$. One can take their direct sum and get a mapping $pr^i_a\oplus pr^i_b\oplus pr^i_c : \mathbb{Z}^{30}\rightarrow (\mathbb{Z}^{6})^{\oplus 3}$. Consider the  mapping

\begin{equation}
\label{pr}
 pr:=\bigoplus_{I=1}^4  pr^i_a\oplus pr^i_b\oplus pr^i_c:   \mathbb{Z}^{30} \oplus   \mathbb{Z}^{30} \oplus   \mathbb{Z}^{30}\oplus   \mathbb{Z}^{30} \rightarrow  \bigoplus_{i=1}^4   (\mathbb{Z}^{6})^{\oplus 3}.
\end{equation}

In the image of $\bigoplus_{i=1}^4(\mathbb{Z}^{6})^{\oplus 3}$ one can introduce a basis $e_{X}^{A^j,i}$, indexed as follows. The lower index is a proper subset of $X\subset \{1,2,3\}$. Thus, these indexes list the basis vectors in the selected $\mathbb{Z}^6$.
There are two upper indexes. The first corresponds to the name of one of the variables (or differential operator), which is substituted in $f_i$ in \eqref{6jsk}. Thus, fixing this index determines the choice of one of the three terms in $(\mathbb{Z}^{6})^{\oplus 3}$. The second upper index is the index $i$, specifying the term $(\mathbb{Z}^{6})^{\oplus 3}$ in $\bigoplus_{i=1}^4(\mathbb{Z}^{6})^{\oplus 3}$.

Thus in the image of  $(\mathbb{Z}^{6})^{\oplus 3}$, corresponding to  $f_1$,  one has the basis 
 
 \begin{align*}
 & \text{In the first summand  $\mathbb{Z}^6$: } e^{A^1,1}_1,  e^{A^1,1}_{2}, e^{A^1,1}_{3} e^{A^1,1}_{1,2}, e^{A^1,1}_{1,3}, e^{A^1,1}_{2,3}.\\
  & \text{In the second summand  $\mathbb{Z}^6$: } e^{A^2,1}_1,  e^{A^2,1}_{2}, e^{A^2,1}_{3} e^{A^2,1}_{1,2}, e^{A^2,1}_{1,3}, e^{A^2,1}_{2,3}.\\
   & \text{ In the third summand   $\mathbb{Z}^6$: } e^{A^4,1}_1,  e^{A^4,1}_{2}, e^{A^4,1}_{3} e^{A^4,1}_{1,2}, e^{A^4,1}_{1,3}, e^{A^1,1}_{2,3}.\\
 \end{align*}
 
In the image of  $(\mathbb{Z}^{6})^{\oplus 3}$, corresponding to  $f_2$,  one has the basis 
 
 \begin{align*}
 & \text{In the first summand  $\mathbb{Z}^6$: } e^{A^4,2}_1,  e^{A^4,2}_{2}, e^{A^4,2}_{3} e^{A^4,2}_{1,2}, e^{A^4,2}_{1,3}, e^{A^4,2}_{2,3}.\\
 & \text{In the second summand  $\mathbb{Z}^6$: } e^{A^3,2}_1,  e^{A^3,2}_{2}, e^{A^3,2}_{3} e^{A^3,2}_{1,2}, e^{A^3,2}_{1,3}, e^{A^3,2}_{2,3}.\\
 & \text{In the third summand  $\mathbb{Z}^6$: } e^{A^5,2}_1,  e^{A^5,2}_{2}, e^{A^4,2}_{3} e^{A^5,2}_{1,2}, e^{A^5,2}_{1,3}, e^{A^5,2}_{2,3}.\\
 \end{align*}
 
In the image of  $(\mathbb{Z}^{6})^{\oplus 3}$, corresponding to $f_3$, one has the basis 

 \begin{align*}
 & \text{In the first summand  $\mathbb{Z}^6$: } e^{A^2,3}_1,  e^{A^2,3}_{2}, e^{A^2,3}_{3} e^{A^2,3}_{1,2}, e^{A^2,3}_{1,3}, e^{A^2,3}_{2,3}.\\
 & \text{In the second summand $\mathbb{Z}^6$: } e^{A^3,3}_1,  e^{A^3,3}_{2}, e^{A^3,3}_{3} e^{A^3,3}_{1,2}, e^{A^3,3}_{1,3}, e^{A^3,3}_{2,3}.\\
 & \text{In the third summand   $\mathbb{Z}^6$: } e^{A^6,3}_1,  e^{A^6,3}_{2}, e^{A^6,3}_{3} e^{A^6,3}_{1,2}, e^{A^6,3}_{1,3}, e^{A^6,3}_{2,3}.\\
 \end{align*}
 
In the image of $(\mathbb{Z}^{6})^{\oplus 3}$, corresponding to  $f_4$,  one has the basis

 \begin{align*}
 & \text{In the first summand  $\mathbb{Z}^6$: } e^{A^1,4}_1,  e^{A^1,4}_{2}, e^{A^1,4}_{3} e^{A^1,4}_{1,2}, e^{A^1,4}_{1,3}, e^{A^1,4}_{2,3}.\\
 & \text{In the second summand  $\mathbb{Z}^6$: }  e^{A^6,4}_1,  e^{A^6,4}_{2}, e^{A^6,4}_{3} e^{A^6,4}_{1,2}, e^{A^6,4}_{1,3}, e^{A^6,4}_{2,3}.\\
 & \text{In the third summand   $\mathbb{Z}^6$: } e^{A^5,4}_1,  e^{A^5,4}_{2}, e^{A^5,4}_{3} e^{A^5,4}_{1,2}, e^{A^5,4}_{1,3}, e^{A^5,4}_{2,3}.\\
 \end{align*}

\begin{defn}
	\label{d} Using the pairing in   \eqref{6jsk}, let us define a sublattice 
 $
 D\subset \bigoplus_{i=1}^4 (\mathbb{Z}^{6})^{\oplus 3}
 $

as a lattice generated vectors for all possible $X\subset\{1,2,3\}$, which are obtained as sums of vectors with the same $X$, the same variable $A^j$, but corresponding to different $f_i$.

 \end{defn}
The lattice   $D$ is generated by vector
 
 \begin{align*}
& e_X^{A^1,1}+e_X^{A^1,4}, \,\,\,     e_X^{A^2,1}+e_X^{A^2,3},   \,\,\,      e_X^{A^3,2}+e_X^{A^3,4}, \\& e_X^{A^4,1}+e_X^{A^4,2} ,   \,\,\,    e_X^{A^5,2}+e_X^{A^5,4},  \,\,\,   
  e_X^{A^6,3}+e_X^{A^6,4} 
 \end{align*}

 Note the following. Let  us be given  a monomial which is obtained from the decomposition of \eqref{6jsk}. It gives a non-zero contribution if for each variable $A^j_X$, $j=1,...,6$ the  exponent of this variable coincides with the exponent of the differentiation by $A^j_X$. The fact of existence of such monomials can be  reformulated in terms of supports. The result is the following.

  \begin{thm}  \label{ot1}

In order for the $6j$-symbol \eqref{6jsk} to be non-zero, it is necessary that

$$
H\cap pr^{-1}( D)\ne \emptyset,
$$
where $H$ is defined in \eqref{h}, and $D$ is defined in the definition of \ref{d}.

  	\end{thm}
 \begin{remark}

This condition can also be formulated in terms of supports of functions $f_i$ as functions of variables $A^j$. It can be seen that the support $supp_{A}f_i$  (a support of a  function of variables $A^j$) is an image of the support $suppf_i$ (a support of a  function of variables $Z$) undert the mapping $pr^i_a\oplus pr^i_b\oplus pr^i_c$, and therefore is also a shifted lattice. Then the necessary condition for $6j$-symbol to be non-zero is formulated as follows:
$
supp_A(f_1f_2f_3f_4)=pr(H)
$
has a non-empty intersection with $D$. To calculate the value of a $6j$-symbol, it is more convenient for us to deal with $Z$ variables.

 \end{remark}
 
 \subsection{A formula for a  $6j$-symbol}
 \label{r33}

 The functions $f_1,f_2,f_3,f_4$ can be considered as functions  of the variables $Z_1,Z_2,Z_3,Z_4$, where $Z_i$ is a set of variables obtained from the set \eqref{pz} by replacing the symbols $a,b,c$ by symbols of variables (or differential operators), substituted in $f_i$ in \eqref{6jsk}. More precise,

\begin{align}
\begin{split}
\label{zmn1}
 & \text{ For  $f_1$: }a,b,c\mapsto A^1,A^2,A^4,\,\,\,\,  \text{ for  $f_2$: }a,b,c\mapsto A^4,A^5,A^6,\\
 & \text{ for  $f_3$: }a,b,c\mapsto A^2,A^3,A^5,\,\,\,\,  \text{ for  $f_4$: }a,b,c\mapsto A^1,A^6,A^5
 \end{split}
\end{align}

The computation of the expression \eqref{6jsk} can be described as follows. First we consider $f_1...f_4$ as  functions of the collection of variables $Z_1,Z_2,Z_3,Z_4$. Let us write $f_1...f_4$ is the sums of the products of monomials in these variables. We leave only the terms whose  support is in $
H\cap pr^{-1}( D)\ne \emptyset.
$ Next, each of the variables from the set $Z_1,Z_2,Z_3,Z_4$, is replaced by the ruler \eqref{zmn} (one takes into account the correspondence of \eqref{zmn1}). Thus we pass to the variables $A^j_X$ or $\frac{\partial}{\partial A_X^j}$ in according to the variable or the differential operator is involved in $f_i$ in \eqref{6jsk}. Multiply them as if they were commuting. After that, in the resulting monomial in $A^j_X$ and $\frac{\partial}{\partial_X^l}$ one applies the differential operators to variables and substitutes zero instead of all the variables.

So, if one focuses on the variable $A_1^1$ that occurs in the monomial, then our actions look as follows. Such a symbol occurs in the notation of variables included in the sets $Z_1$ and $Z_4$. We take the monomial resulting from the decomposition of $f_1...f_4$. Let its support belong to $
H\cap pr^{-1}( D)$. Let's write it explicitly together with the coefficient at this monomial. It was noted above that $f_1$,...,$f_4$ are $\Gamma$  series in the  variables $Z_1,Z_2,Z_3,Z_4$, so that the coefficient at the monomial is a product of the inverse quantities to the factorials of degrees:

\begin{equation}
\label{monom}
\underbrace{\frac{[A^1_1A^2_{2,3}]^{\beta_1}}{\beta_1!}     \frac{[A^1_1A^4_{2,3}]^{\beta_2}}{\beta_2!}...         }_{\text{РёР· }f_1}\cdot...\cdot
\underbrace{\frac{[A^1_1A^6_{2,3}]^{\gamma_1}}{\gamma_1!}     \frac{[A^1_1A^5_{2,3}]^{\gamma_2}}{\gamma_2!}...         }_{\text{РёР· }f_4}.
\end{equation}

Next,  one calculates the sum of the  exponents of  the variables, whose notation contains the symbol $A_1^1$. For factors originating from $f_1$, this sum is equal to $\beta_1+\beta_2+..$, and for factors originating from $f_4$, this sum is equal to $\gamma_1+\gamma_2+...$. The fact that the supports belong to $
H\cap pr^{-1}(D)$ implies that $\beta_1+\beta_2+...=\gamma_1+\gamma_2+...$.
When passing to $A^j_X$ or $\frac{\partial}{\partial A_X^j}$, one substitutes  $\frac{\partial}{\partial A_1^1}$ into factors  originating from $f_1$, and  one substitutes   $A_1^1$  in factors  originating from $f_4$. After applying the differential operator to the variable $A_1^1$ and substituting zero instead of $A_1^1$ in \eqref{monom}, all symbols $A_1^1$ are actually deleted and a numerical multiplier $(\beta_1+\beta_2+...)!$ is added in the top.

After performing similar actions with all  the variables $A_X^j$, the monomial \eqref{monom} turns into a numeric fraction. In its denominator occur the factorials (in the multi-index sense) of the degree of the monomial as a monomial of the variables $Z_1,...,Z_4$ and at the top occur the factorials of the degree (again in the multi-index sense)  of the monomial as a monomial in the variables $A_X^j$.

Let us proceed to the formulation of the main result of this paper. The set $
H\cap pr^{-1}(D)$ is a shifted lattice. Therefore, for some vector $\varkappa$ and some lattice $L\subset\mathbb{Z}^{30}\oplus\mathbb{Z}^{30}\oplus\mathbb{Z}^{30}\oplus\mathbb{Z}^{30}$, one can write

\begin{equation}
\label{kp}
H\cap pr^{-1}( D)=\varkappa+L\subset  \mathbb{Z}^{30}\oplus \mathbb{Z}^{30}\oplus\mathbb{Z}^{30}\oplus\mathbb{Z}^{30}
\end{equation}

There is a projection  $pr$ defined by the formula \eqref{pr}. One associates with the shifted lattice $\varkappa+L$ a series of hypergeometric type (in fact, being a finite sum) in the variables $Z=\{Z_1,...,Z_4\}$, defined by a formula related to \eqref{gmr}:

\begin{equation}
\mathcal{J}_{\gamma}(Z;L)=\sum_{x\in \varkappa+L} \frac{\sqrt{\Gamma(pr(x)+1)}Z^{x}}{\Gamma(x+1)}
\end{equation}

\begin{thm}\label{ot2}
The $6j$-symbol \eqref{6jsk} is equal to $\mathcal{J}_{\gamma}(\pm 1;L)$, where $1$ is substituted instead of those variables that are entering  in the determinants in $f_i$ the sign $+$ and $-1$ is substituted instead of those variables that are entering in the determinants with the sign $-$.

The shifted lattice  $ \varkappa+L$ is defined in  \eqref{kp}, where  $H=\oplus H_i$, $H_i=suppf_i$; $D$ is given in the Definition  \ref{d}, and  $pr$ is defined in \eqref{pr}.

\end{thm}

\end{document}